\UseRawInputEncoding
\documentclass[reqno,12pt,a4paper]{amsart}
\usepackage{tgbonum}
\usepackage[margin=2.5cm]{geometry}

\usepackage{hyperref}
\usepackage{longtable,fullpage}
\usepackage{booktabs}
\usepackage{colortbl}
\usepackage{arydshln}
\usepackage{verbatim}
\usepackage{multirow}
\usepackage{makecell}
\usepackage{pdfsync}

\usepackage{amsthm}
\usepackage{hyperref}
\usepackage{enumerate}
\usepackage{url}
\usepackage{pdflscape}
\usepackage{longtable}
\usepackage{booktabs}
\usepackage{colortbl}
\newcommand{\evnrow}{\rowcolor[gray]{0.95}}

\usepackage[margin=2.5cm]{geometry}

\usepackage{pdflscape}
\usepackage{longtable}
\usepackage{booktabs}
\usepackage{colortbl}
\usepackage{arydshln}
\usepackage{verbatim}
\usepackage{enumitem}
\usepackage{mathdots}

\usepackage[english]{babel}
\usepackage{lipsum}
\usepackage{amsthm}
\usepackage{todonotes}
\usepackage[vcentermath]{youngtab}

\theoremstyle{plain}
\newtheorem{lema}{Lemma}[section]

\newtheorem{thrm}[lema]{Theorem}

\theoremstyle{remark}

\newtheorem{rmk}[lema]{Remark}

\theoremstyle{definition}
\newtheorem{dfn}[lema]{Definition}
\newtheorem{ej}[lema]{Example}

\def\ZZ{{\mathbb Z}}
\def\QQ{{\mathbb Q}}
\def\CC{{\mathbb C}}

\def\AA{{\mathbb A}}

\def\PP{{\mathbb P}}
\def\wgr{\mathrm{wGr(2,5)}}
\def\Gr{\mathrm{Gr}}
\def\wp2{\mathrm{w}(\PP^2\times \PP^2)}

\def\wG{\mathrm{w}\mathcal{G}}
\def\wP{\mathrm{w}\mathcal{P}}

\def\cF{\mathcal{F}}

\def\sm{\mathrm{smooth}}

 \newcommand{\into}{\hookrightarrow}

 \newcommand{\Oh}{\mathcal O}

 \newcommand{\sB}{\mathcal B}

 \newcommand{\Si}{\Sigma}

\newcommand{\PxP}{\PP^2 \times \PP^2}
\usepackage{amstext}
\usepackage{amssymb}

\begin{document}

%--------------------------------------------------------------------------------------------

\author[M.I.~Qureshi]{Muhammad Imran Qureshi}
\address{Deptartment of Mathematics, King Fahd University of Petroleum and Minerals, Dhahran, Saudi Arabia, 31261 
}
\address{ Interdisciplinary Research Center for Intelligent Secure Systems,  King Fahd University of Petroleum \& Minerals, Dhahran, Saudi Arabia, 31261} 
\email{imran.qureshi@kfupm.edu.sa, i.qureshi@maths.oxon.org}
%-------------------------------------------------------------------------------
%\subjclass[2010]{14J30 (Primary); 14Q15, 14M07 (Secondary)}
%--------------------------------------------------------------------------------------------
\title{Terminal  Fano four folds in low codimension}
%--------------------------------------------------------------------------------------------
\begin{abstract}
We construct well-formed and quasismooth terminal Fano 4-folds of index 1 in low codimension containing at worst isolated orbifold points. We provide a certain classification of these varieties  where their images under the anitcanonical embedding can be described as  codimension 2, 3, or 4 subvarieties of some weighted projective space. In particular, we focus on isolated terminal Fano 4-folds that either have an empty linear system or a relatively large one, but whose linear section is not an isolated canonical Calabi--Yau 3-fold. In total, we classify 95 families of terminal Fano 4-folds of the first type and 32 families of the second type. We also describe our algorithmic approach and the pivotal role of computer algebra in our results.

\end{abstract}
%--------------------------------------------------------------------------------------------
\maketitle
%\tableofcontents
%--------------------------------------------------------------------------------------------
\section{Introduction}
\subsection{ Terminal Fano varieties} A  \textit{Fano variety}  \(X\) is an algebraic variety with an ample anitcanonical divisor class \(-K_X\). 
The \textit{Fano index} \( i \) is defined as the largest integer \( i \) such that the anitcanonical divisor class \( -K_X \) can be written as \( -K_X = i \cdot D \) for some $\QQ$-ample Weil divisor \( D \) on \( X \). A \textit{quotient singularity}
of type  $\frac{1}{r}(a_1, \ldots, a_n)$ is a quotient of \(\AA^n\) by a
cyclic group  \(\mu_r\) of order \(r\) given \[\epsilon: x_i \mapsto e^{a_i}x_i,
\; i=1,\ldots,n, \] where $r$ is a positive integer and $a_1, \ldots, a_n \in \mathbb{Z}_{\ge 0}$. It is called \textit{isolated} if  all \(a_i\)s  are coprime to $r$ and it is called \emph{terminal} if and only if 
\[
\frac 1r\sum_{i=1}^n {\overline{ka_i}}>1 \text{ for all } k=1,\ldots,r-1
.\]
 If \(X\)  contains at worst  terminal singularities then
it is called a terminal Fano variety.

The classification of Fano varieties is a central focus in studying algebraic varieties. It is well established that there are only finitely many deformation families of smooth Fano varieties in each dimension \cite{KMM}. The classification has been completed for dimensions less than or equal to three \cite{Isk1, Isk2, Isk-Pro, MM}.  
For dimensions greater than or equal to 4, the comprehensive classification remains unfinished. The complete classification of smooth Fano 4-folds with an index of 2 or greater is established, and there are 35 deformation families of such Fano 4-folds \cite{Fuj1,Fuj4,Isk1,Isk-Pro,Kob1,Wilson,Wis,Q-periods}. In index  one case, there are partial classification
of results  \cite{Kuchle-95,Kuchle-97,Bat-toric-4,CKP,Klash,HLM,QBAM}, leading
to around one thousand  smooth families of Fano 4-folds. 
 
In dimension greater than or equal to three,   terminal singularities \cite{Reid-MM}, are essential to construct minimal models of algebraic varieties. Therefore,
the classification of terminal Fano varieties has gained much traction over
the years.  The graded ring database \cite{GRDB,GRDB-paper} lists all possible 
Hilbert series and further expected attributes (terminal singular points, equation degrees,
ambient weighted projective space, Hilbert numerator, etc.) that may correspond to  a terminal \(\QQ\)-Fano 3-fold.  A lot of
progress has been made in classifying these varieties in
 \cite{c3f,fletcher,ABR,Takagi,BKR,BKQ,CD-Fano}.   
In dimension four,  terminal Fano 4-fold weighted projective spaces have been classified by Kasprzyk in \cite{Al-twps4}.
The case of  terminal Fano 4-fold hypersurfaces of index 1 is classified by Brown and Kasprzyk \cite{BK-ths4} and there are  a total of $490$ families of isolated terminal Fano 4-folds (ITF4s) on the graded ring database \cite{GRDB}. 
%%%%%%%%%%%%%%%%%%%%%%%%%%%%%%%%%%%%%%%%%%%%%%%%%%%%%%%%%%%%%%%

\subsection{Main Results} In this article, we  describe the construction and classification  of ITF4s embedded in some weighted projective
space as an algebraic variety of  codimension greater than or equal to 2. We study families of  ITF4s such that   under their anitcanonical   embeddings they can be described  as wellformed and quasismooth varieties of   codimensions 2, 3, or 4 in some weighted projective space.
In particular, we focus on two types of ITF4s in this article, for which we provide our motivation in Section \ref{motivation}, which we denote by type-$K_0$ and type-$K_2$.  

\subsubsection{Type-$K_0$} 
We call a class of  ITF4s of index 1  to be of \emph{type-$K_0$} if its   linear  system $|-K_X|$ is empty, i.e.  \(\left( h^0(-K_X)=0\right)\). Consequently,  it does not have any   Calabi--Yau 3-fold as its hyperplane  section.  

\begin{thrm}[Classification of type-$K_0$ isolated terminal Fano fourfolds]\label{Th:type-I}
Let \(X\) be a {wellformed} and {quasismooth} isolated terminal Fano 4-fold (ITF4) of  index \(1\)  with $h^0(-K_X)=0$.
If \(X\) admits an anticanonical embedding
\[
X \;\hookrightarrow\;
\PP^{\,4+c}\!\bigl(w_{0},\dots,w_{4+c}\bigr),
\qquad c\in\{2,3,4\},
\]
that can be described by a Gorenstein format \(\mathcal{F}\) of codimension \(c\), satisfying \(\sum_{i=0}^{4+c} w_{i}\;\le\; W\).  Then  in each Gorenstein format $\cF$ there exist at least  \(m\)  families of such ITF4s, summarized in Table~\ref{tab-summary-t1}.
\end{thrm}

\begin{center}
\begin{table}[h]
\caption{ITF4s of  Type-$K_0$}   
\label{tab-summary-t1}
\renewcommand*{\arraystretch}{1.5}
\begin{tabular}{ |c|c|c|c|c| }
  \hline
 \multirow{2}{*}{Attributes} & \multicolumn{4}{|c|}{Format \(\mathcal F\)} \\\cline{2-5}
  & C.I & C.I & Gr(2,5) & $\PxP$  \\
   \hline\hline 
  \evnrow \(c\) & 2 & 3 & 3 & 4 \\
   \(m\) & 80 & 13 & 1 & 1 \\ 
  \evnrow \(W\) & 101 & 70 & 71 & 59 \\
  \hline
\end{tabular}
\end{table}
\end{center}

\subsubsection{Type-$K_2$} We call a class of  ITF4s of index 1 to be of \emph{type-$K_2$} if the linear system  $|-K_X|$ satisfies   \(h^0(-K_X)\ge 2\), yet a general member \(Y\in|-K_{X}|\) fails to be a Calabi–Yau threefold with canonical isolated
orbifold points. In our case, we consider those  ITF4s  $X$ that contain at least  one orbifold point of type \[
  \tfrac{1}{r}\bigl(a_{1},a_{2},a_{3},a_{4}\bigr), \qquad a_{i}\neq1\;(i=1,\dots,4).
\]
Therefore, these families are  new and do not occur in the list of canonical Calabi–Yau 3-fold (cCY3)
sections of Fano 4-folds  in~\cite{BKZ}.

\begin{thrm}[Classification of type-$K_2$ isolated terminal Fano fourfolds]\label{Th:type-II}
Let \(X\) be a wellformed, quasismooth isolated terminal Fano 4-fold (ITF4) of index~\(1\)
with  \(h^{0}(-K_{X})\ge 2\) and whose
general anticanonical  section is not  an isolated canonical Calabi–Yau 3-fold. Assume \(X\) admits an anticanonical embedding
\[
X \;\hookrightarrow\;
\PP^{\,4+c}\!\bigl(w_{0},\dots,w_{4+c}\bigr),
\qquad c\in\{2,3,4\},
\]
that can be described by a Gorenstein format \(\mathcal{F}\) of codimension \(c\), satisfying \(\sum_{i=0}^{4+c} w_{i}\;\le\; W\).  Then  in each Gorenstein format $\cF$ there exist at least  \(m\)  families of such ITF4s, summarized in Table~\ref{tab-summary-t2}.
\end{thrm}
\begin{center}
\begin{table}[h]
\caption{ITF4s of Type-$K_2$}   
\label{tab-summary-t2}
\renewcommand*{\arraystretch}{1.5}
\begin{tabular}{ |c|c|c|c|c| }
  \hline
 \multirow{2}{*}{Attributes} & \multicolumn{4}{|c|}{Format \(\mathcal F\)} \\\cline{2-5}
  & Com. Int. & Com. Int. & Gr(2,5) & $\PxP$  \\
   \hline\hline 
   \evnrow\(c\) & 2 & 3 & 3 & 4 \\
   \(m\) & 12 & 8 \ & 7 & 5 \\ 
   \evnrow\(W\) & 101 & 70 & 71 & 57 \\
  \hline
\end{tabular}
\end{table}
\end{center} 
\begin{rmk}
 In other words, if 
 \(X\subset\PP^{\,4+c}\bigl(w_{0},\dots,w_{4+c}\bigr)\) 
with \(\sum w_{i}\le W\) be an index-\(1\) isolated terminal Fano 4-fold that satisfies the
assumptions of Theorem~\ref{Th:type-I} or Theorem~\ref{Th:type-II} in one of the prescribed
Gorenstein formats.
Then \(X\) is isomorphic to exactly one of the \(95\) Type-$K_0$ or \(32\) Type-$K_2$ families recorded in
Tables~\ref{tab-summary-t1} and~\ref{tab-summary-t2}
(the complete data are also available on GitHub\footnote{\url{https://github.com/QureshiMI/TerminalFano4}}).
Conversely, each of the \(95+32\) cases listed in those tables is realized by an ITF4
with the stated properties in the corresponding Gorenstein format.
\end{rmk}

\begin{rmk}
 We also searched for examples of ITF4s as codimension-4 weighted complete intersections, but no quasismooth example was found.  We searched up to \(W = 64\) in that case and found 13 candidate isolated terminal Fano 4-folds, with seven of those satisfying \(h^{0}(-K_{X}) \ge 2\) and none   has an empty linear system.  Except for the well-known smooth Fano 4-fold complete intersection \(X_{2,2,2,2} \subset \mathbb{P}^{8}\),  all of the candidates have non-terminal orbifold points. 

\end{rmk}
  We use a format search algorithm of \cite{QJSC,BKZ}, given in
Section \ref{algorithm}, to search for the candidate ITF4s.  Then we prove the existence of each 4-folds by using a combination of theoretical and computational tools, explained with details in Section\ref{Flowchart}.
The above list is indeed not a complete classification of ITF4s which can be constructed in these formats and  undoubtedly a sublist of all possible cases, however, it is a complete list up to the given value of \(W\).
The  format search   parameter, i.e. the sum of the weights \(W\) on the weighted
projective space containing a terminal Fano 4-fold, is unbounded. Thus,  we search for examples for up to certain
value of \(W\), determined by the exhaustion of  computer memory for the given format.
\subsection{Context and Motivation}
\label{motivation}
 One way to look at our results as a natural extension of already existing
classification results  of terminal Fano 3-folds and 4-folds. However, the main motivation to study  these two classes of ITF4s is to compare these with the  existing geography of 3-dimensional case in codimension \(\le 4\), available on the graded ring database \cite{GRDB-paper}.  We make the comparison   for each type of ITF4s separately.

If we compare the type-$K_0$ case then in codimension \(\le 4\), there are only four examples of terminal \(\QQ\)-Fano 3-folds with the empty linear system. One of these is a codimension 2 weighted complete intersection and 3 of them can be described as regular pullbacks in some cluster Gorenstein format \cite{CD-Fano}. But  in dimension four, we notice that  there is a large set of such terminal Fano 4-folds in each codimension. In total, we obtain 95 examples of such 4-folds in this article. \ There is no example in 3-fold case with \(h^{0}(-2K_X)=0\) but in dimension 4 we have an  example \ref{extreme-CI2} of codimension 2 terminal Fano 4-fold  with where first four plurigenera are zero, i.e.   \(h^0(-lK_X)=0,\text{ for } 1\le l\le 4.  \)  We provide the detailed results
of type-$K_0$ for each format in Table \ref{tab-summary-KX-empty}. 
% \begin{center}
% \begin{table}[h]
% \caption{ Number of ITF4s with \(h^{0}(-lK_X)=0\) for \( l< M\)   }   
% \label{tab-summary-H0}
% \renewcommand*{\arraystretch}{1.5}
% \begin{tabular}{ |l|l|l|l|l|l|l|l|l |l}
%   \hline
%  \multirow{2}{*}{$M$} &\multicolumn{9}{|c|}{Format \(\mathcal F\)} \\\cline{2-5}
%   & C.I cod. 2 &C.I cod. 3&Gr(2,5)&$\PxP$  \\
%    \hline\hline 
% \evnrow   \(2\) &61&6&1&1\\
%    3 &16&1 &0&0\\ 
% \evnrow   4 &2&1&0&0\\ 
%    5 &1&0 &0&0\\ 
%    \evnrow 
% 
%    \hline
% 
% \end{tabular}
% \end{table}
% 
% \end{center}
% 
                
For the comparison in the type-$K_2$ case, recall that   a terminal \(\QQ\)-Fano 3-fold \(V\) of index 1  has quotient singularities of the form \(\frac1r(1,a,b)\) where \(a,b\) are relatively prime to \(r\). If \(|-K_V|\ne \emptyset  \) then its linear section is always a  \(K3\) surface (a two Calabi--Yau) with Du Val singularities \(\frac1r(a,b)\), which are two-dimensional canonical quotient singularities. However, we observe that in the case of   terminal \(\QQ\)-Fano 4-folds this doesn't happen necessarily, i.e there exist many families of ITF4s with \(|-K_X|\ne \emptyset\) but a generic linear section is not a canonical Calabi--Yau 3-fold (cCY3) with isolated orbrifold points. This phenomenon is relatively widespread and we  focus on such cases with \(h^0(-K_X)\ge 2.  \) In total
we obtain at least 32 such families of ITF4s and we provide a detailed division
in various formats in Table \ref{Combined-Summary-Fano4-Sorted}.

  %-------------------------------------------------------------------------
% 
% \begin{table}[h]
% \caption{Number of ITF4s with \(h^{0}(-K_X) \ge 2\) and linear section is not a cCY3}
% \label{tab-summary-NCY3}
% \renewcommand*{\arraystretch}{1.5}
% \centering
% \begin{tabular}{|l|c|c|c|c|}
%   \hline
%   \multirow{2}{*}{\(h^{0}(-K_X)\)} & \multicolumn{4}{c|}{Format \(\mathcal{F}\)} \\ \cline{2-5}
%    & C.I codim. 2 & C.I codim. 3 & Gr(2,5) & \(\PxP\) \\ 
%   \hline\hline
%   2 & 5 & 5 & 5 & 5 \\ 
%   \hline
%   3 & 0 & 1 & 1 & 3 \\ 
%   \hline
%   4 & 0 & 0 & 0 & 1 \\ 
%   \hline
% \end{tabular}
% \end{table}

\section{Preliminaries} 
In this section, we provide some basic definitions and the necessary background
required.
 
\noindent A weighted projective variety  \(X \into \PP^N(w_i)\) of codimension \(c\) is called \textit{wellformed} if it does not contain a singular stratum of   codimension \(c+1\). 
 It is called \textit{quasismooth} if the affine cone over \(\tilde X\) in \(\AA^{n+1}\)
is
 smooth outside the vertex of the ambient affine space. 
 A \textit{format} is roughly a way of representing the equations of varieties systematically. For example, the Segre embedding of \(\PxP\)  can be described as \(2 \times 2 \) minors of the size 3 matrix. A more formal definition of Gorenstein format  is given below. 
\begin{dfn}\cite{BKZ}\label{dfn:GF} A codimension \(c\) \emph{Gorenstein format}  \(\mathcal F\) is a triple \(\left(\widetilde{V } 
 , \mathcal{R}, \mu \right)\) which consists of a codimension \(c\) affine Gorenstein variety \(\widetilde V\subset \AA^n\), a minimal graded free resolution \(\mathcal R\) of \(\Oh_{\widetilde V}\) as a graded \(\Oh_{\AA^n}\) module, and a \(\CC^*\)-action \(\mu\) of strictly positive weights on \(\widetilde V \).  
 \end{dfn}
 We only consider those Gorenstein formats where the action \(\mu\) leaves the variety \(\widetilde V\) invariant  and the  free resolution \(\mathcal R\) is equivariant for the action. The varieties defined below in    \ref{dfn:WG} and \ref{dfn:P2xP2} are examples of such Gorenstein formats.
\begin{dfn}\cite{fletcher} A weighted projective variety \(X_{d_1,d_2,\ldots,d_c}\into\PP(w_0,\ldots,w_n) \) defined by \(c\) weighted homogeneous polynomials of degrees \(d_1,\ldots,d_c\) is called \textit{weighted complete intersection} of codimension \(c\) if its dimension is equal to \(n-c\). If  \(d_j\ne w_i\) for all  \( 1\le j\le c\) and \(0\le i\le  n\), then
we say that is not the intersection with a linear cone of \(\PP(w_i)\). We
only consider such weighted complete intersections in this paper.
\end{dfn}
\begin{dfn}\cite{wg}\label{dfn:WG} The Grassmannians of 2-planes in \(\CC^5\)
has the  Pl\"ucker embedding 
  \(\Gr(2,5) \into\PP^9\left(\bigwedge^2 \CC^5\right) \). To define a weighted
  Grassmannian,   we choose a vector of half integers   $w:=(c_{1},\cdots,c_5)$   satisfying   $$c_i+c_j>0,\; 1\le i<j\le 5.$$ Then the quotient of affine cone of the Pl\"ucker embedding  minus the vertex     $\widetilde{\Gr(2,5)} \backslash\{\underline 0\} $ by $\CC^\times$  given by: 
   $$\mu:x_{ij}\mapsto \mu^{c_i +c_j }x_{ij}$$ \end{dfn}
   defines the \textit{weighted Grassmannian} \(\wgr\), embedded in the weighted
   projective space 
\begin{equation}\label{embed} \PP(\left\{w_{ij}: 1\le i <j\le 5, w_{ij}=c_i+c_j\right\}).  \end{equation} We may use \(\wG\) to denote the \(\wgr\). The image of \(\Gr(2,5)\) (and indeed of \(\wG\) ) inside \eqref{embed} can be defined by  five  $4\times 4$  Pfaffians of $5\times 5$ skew symmetric matrix  
   \begin{equation}\label{eq:pf_mat}\left(\begin{matrix} 
x_{12}&x_{13}&x_{14}&x_{15}\\ 
&x_{23}&x_{24}&x_{25}\\ 
&&x_{34}&x_{35}\\ 
&&& x_{45} \end{matrix}\right)
,\end{equation}
where we omit the the lower triangular part and diagonal of zeros  of the skew symmetric matrix.  
If   $\wG$ is wellformed,   then the  canonical divisor  class is given by  \begin{equation}\label{eq:K_G}K_{\wG}= \left(-\frac12\sum_{1\le i<j\le 5} w_{ij} \right)H,\end{equation}
for an ample divisor \(H\). 
  
\begin{dfn}\cite{BKQ} \label{dfn:P2xP2} Let \(\Si\) be the Segre embedding of  \begin{equation}\PxP\into \PP^8(x_{ij}),\;\; 1\le i,j\le 3\end{equation}  then the equations of the image can be
described by \(2\times 2 \) minors of  \(3\times 3\) matrix \(x_{ij}\).
We take a pair of  half-integer 3-tuples   $a=(a_1,a_2,a_3)$ and $b=(b_1,b_2,b_3)$
such that  
$$a_i+b_j > 0,a_i\le a_j \text{ and } b_i\le b_j \text{ for } 1\le i\le j \le 3. $$  
Then   the we define the\textit{ weighted} \( \PxP\) to be the quotient  :
$$\mu:x_{ij}\mapsto\mu^{a_i+b_j}x_{ij}, \; 1\le i,j\le 3,$$ of the punctured affine cone $\widetilde{\Si}\backslash\{\underline 0\}$  by  $\CC^\times$.
We may also  denote by \(\wP\). Thus for a choice of parameter $p=(a_1,a_2,a_3;b_1,b_2,b_3), $ we get the embedding $$\wP\into \PP^8\left(w_{ij}: w_{ij}=a_i+b_j;1\le i,j\le 3\right). $$ 
\end{dfn}
We usually use the   weight matrix 
\begin{equation}\label{eq:wtmx_P2}
\begin{pmatrix} w_{11} & w_{12} & w_{13} \\ w_{21} & w_{22} & w_{23} \\ w_{31} & w_{32} & w_{33} \end{pmatrix} \text{ where }w_{ij}=a_i+b_j; 1\le i,j\le 3,
\end{equation}
to refer to a specific \(\wP\). 
 If $\wP$ is wellformed then the  canonical divisor class is given by \begin{equation}\label{eq:K_P2}K_{\wP}=\left(-\sum_{ i=1}^3 w_{ii} \right)H,\end{equation}
for an ample divisor \(H\).
%%%%%%%%%%%%%%%%%%%%%%%%%%%%%%%%%%%%%%%%%%%%%%%%%%%%
\section{Flowchart of the proof}
\label{Flowchart}
In this section, we describe all the ingredients required to prove the Theorem
\ref{Th:type-I} and Theorem \ref{Th:type-II}. 
We use the format search algorithm of \cite{QJSC,BKZ} to search for the candidate ITF4s.  Then we prove the existence of each 4-folds by using a combination of theoretical and computational tools.
\subsection{Algorithm for candidate isolated orbifolds}
\label{algorithm}
This section reviews the algorithm from \cite{QJSC}, which we utilized to generate the list of candidate terminal Fano 4-folds. A key component is the orbifold Riemann--Roch formula of Buckley--Reid--Zhou\ \cite{BRZ}, which provides a formula to decompose the  Hilbert series \(P_X(t)\) of \(X\) as a sum of a smooth component and an orbifold component. It states that if an orbifold \(X\) contains a set \(\mathcal B = \{k_i \times Q_i: m_i \in \ZZ_{>0}\}\) of isolated orbifold points, then:
\begin{equation}\label{hs-decomp}
P_X(t) = P_\sm(t) + \sum_{Q_i \in \mathcal B} k_i P_{Q_i}(t),
\end{equation}
where \(P_{\sm}\) denotes the smooth part, and $\sum_{Q_i \in \mathcal B} k_i P_{Q_i}(t)$  represents the orbifold component of the Hilbert Series. 

 The algorithm generates a comprehensive list of orbifolds with a specified Fano index, dimension, and canonical class \(K_{X} = \Oh(kD)\)  for
a fixed value of adjunction number, determined by the Hilbert series. We  outline the algorithm steps for Fano 4-folds below.

\begin{enumerate}
\item Determine the Hilbert series and canonical class of the ambient
weighted projective variety.
\item Identify all possible embeddings in \(\PP^8(w_i)\) with index \(i\), i.e., \(K_X = \Oh(-i)\), by enumerating weights and utilizing the adjunction formula.
\item For each embedding from step (ii), calculate the Hilbert series of \(X\) and determine the term representing its smooth component \(P_\sm(t)\).
\item Compile a list of all potential terminal  orbifold points from the weights of \(\PP^8(w_i)\) that could be on \(X\).
\item Form all potential subsets of orbifold points, and for each set \(\mathcal B\), assess whether the difference \(P_X(t) - P_{\sm}(t)\) can equate to \(\sum_{Q_i \in \mathcal B} k_i P_{Q_i}(t)\) for some \(k_i\in \ZZ\).
\item \(X\) qualifies as a candidate terminal Fano 4-fold with a basket of singularities \(\mathcal B\) if each  \(k_i \ge 0\). Repeat steps (iii) to (vi) for each embedding determined in step (ii).
\end{enumerate}

\label{steps-proof}
  
\subsection{Existence, consistency and wellformedness} 
Here we provide a brief overview from \cite{QMOC} where more detailed explanation can be found. The  existence of  each candidate Fano 4-fold \(X\) of index 1  as an algebraic variety follows as the equations are obtained from the    equations of the
given Gorenstein
format. Then we check for consistency of the candidate Fano 4-fold with the given  Hilbert series decomposition; that is, we verify that $X$ contains the exact set of terminal singularities predicted by algorithm  \ref{algorithm}. This requires
  computing the intersection of  \(X\) with each singular toric stratum of
 the ambient weighted projective space \(\PP^n(w_i)\), by using the computer algebra system \textsc{MAGMA}.    Since we allow at worst isolated orbifold points, the dimension of each intersection must be less than or equal to zero. Thus, consistency also implies that $X$ is well-formed.  A candidate 4-fold  fails  if either there is an isolated orbifold that is   non-terminal or   it  contains some higher dimensional orbifold locus of \(\PP^n(w_i)\).
 
\subsection{Proving Quasismoothness}  The key  and relatively harder step is to  prove the  quasismoothness  of the candidate Fano 4-fold.  Since 
taking complete intersection with a form of degree \(d\) induces a base locus
of the linear system \(|\Oh(d)|\) that can be quite large due to the intricate
combinatorics of weights and equation degrees of \(X\). By a version of Bertini's theorem, a general member is quasismooth outside of the reduced part of the base locus.  If the base locus happens to be zero dimensional then we can use theoretical arguments to prove quasismoothness. In most cases, the base locus is complicated and  higher dimensional, so we use the  computer algebra system $\textsc{Magma}$ \cite{magma} to prove the quasismoothness of \(X\) by writing   down sparse
representations of its explicit
equations and then apply the Jacobian criteron.

\section{Isolated terminal Fano 4-folds with \(|-K_X|\) empty}

This section consists of the sample examples of ITF4s obtained from  Theorem
\ref{Th:type-I}.  We  present examples that represent the extreme cases, i.e. examples with maximum \(m\) such that \(h^0(-lK_X)=0\) for \(1\le l\le m\). We first provide a tabular summary of  the computer search and examples obtained in each case in Table \ref{tab-summary-KX-empty}.     
 
 \label{Fano4-K-empty}
 %%%%%%%%%%%%%%%%%%%%%%%%%%%%%% 
 \subsection{Table summary }
In this section, we provide the details of the distribution of the terminal Fano
4-folds with the empty linear system in each Gorenstein format, in Table \ref{tab-summary-KX-empty}. The following explains the notations used in the table.

The column labelled "\#Candidates" lists, for each format, the total number of isolated
terminal Fano $4$-fold candidates produced by the algorithm and subsequently analysed.
The maximum value of the sum of the weights  in that format is denoted by~$W$.
The next four rows,
labelled "\#\,{$h^{0}(-\!lK_{X}) = 0,\; l \le n$}",
record how many of those candidates satisfy
$h^{0}(-\!lK_{X}) = 0$ for every $l\le n$.
Finally, the row "QS~Examples" gives the number of candidates that are verified to be
quasismooth.

\begin{table}[h]
\caption{Summary of Isolated Terminal  Fano 4-folds with \( |-K_X| \) Empty}
\label{tab-summary-KX-empty}
\renewcommand*{\arraystretch}{1.5}
\centering
\begin{tabular}{|l|c|c|c|c|}
  \hline
  \multirow{2}{*}{{Attribute}} & \multicolumn{4}{c|}{{Format}} \\ \cline{2-5}
   & {C.I cod. 2} & {C.I cod. 3} & {Gr(2,5)} & \textbf{\( \PxP \)} \\
  \hline\hline

\evnrow  \#Candidates & 702 & 78 & 295 & 176 \\
  
  W & 101 & 70 & 70 & 57 \\
  
\evnrow  \( h^0(-K_X) = 0 \) & 61 & 7 & 1 & 1 \\
  
  \( h^0(-lK_X) = 0,\; l \le 2 \) & 16 & 3 & - & - \\
  
  \evnrow \( h^0(-lK_X) = 0,\; l \le 3 \) & 2 & 3 & - & - \\
  
  \( h^0(-lK_X) = 0,\; l \le 4 \) & 1 & 0 & - & - \\
  
  \evnrow QS Examples  & $80$ & $13$ & 1 & 1 \\
  \hline

\end{tabular}
\end{table}

\subsection{Sample Examples}
In each Gorenstein format, we list an example that represents the extreme case, i.e., the examples with largest \(m\) with \(h^0(-lK_X)=0, 1\le l<m\).
\begin{ej}\label{extreme-CI2}
In the codimension 2 complete intersection format, the extreme case appears for \(m=5\). The weighted complete intersection \[X_{36,40}\into \PP(5^2,7,8,9,12,31)\] is a wellformed and quasismooth  ITF4 with the following basket of orbifold points. 
$$\sB=\{1/3(1, 2, 2, 2), 8 \times1/5(2,2 , 3, 4), 1/7(2, 3, 5, 5), 1/31( 5, 7, 
8,12)\}.$$ 
\end{ej}
\begin{ej}
In the codimension 3 complete intersection format, the extreme case appears for \(m=4\). The weighted complete intersection \[X_{16,18,20}\into \PP(4,5^2,6,7,8,9,11)\] is a wellformed and quasismooth  ITF4 with the following basket of orbifold points. 
$$\sB=\{1/3(1, 2, 2, 2), 4 \times1/5(1,2,4, 4), 1/7(1,4, 5, 5), 1/11( 4,5, 6, 
8)\}.$$ 
\end{ej}
The next example is an ITF4 with an empty linear system in $\Gr(2,5)$ format. 
\begin{ej}
Let \[X_{ 36, 32, 30, 27, 17}\into \PP(2, 3, 5, 7, 8, 9, 11, 27)\] be a candidate Fano
4-fold with weight matrix  
$$\left(\begin{matrix} 3&5&8&18\\ &9&12&22\\ &&14&24\\&&&27 \end{matrix}\right).$$ Then we prove that $X$ is a wellformed and quasismooth terminal Fano 4-fold with the following basket of singularities.

\[
\mathcal{B}
  = \left\{
      \begin{array}{c}
        3\times\tfrac{1}{3}(1,2,2,2),\;
        \tfrac{1}{5}(1,3,3,4),\;
        \tfrac{1}{7}(2,2,5,6), \\[4pt]
        \tfrac{1}{11}(2,5,7,9),\;
        \tfrac{1}{27}(2,7,8,11)
      \end{array}
    \right\}.
\]
\end{ej}
The next example is an ITF4 with an empty linear system in $\PxP$ format. 
\begin{ej}
 Let $X\into \PP(2,3^3,4,5^2,7,11)$ be a Fano 4-fold as regular pullback in \(\PxP\) format given by the weight matrix 
$$\left(\begin{matrix} 3&4&5\\ 5&6&7\\ 10&11&12 \end{matrix}\right).$$
Then the algorithm computesthe basket of singularities to be $$\sB=\{14 \times1/3(1, 2, 2, 2), 1/5(2,3, 3, 3), 1/7(2, 3,5, 5), 1/11(2,3,3,4)\}.$$
\end{ej}

Let
\[
R \;=\;
k\!\bigl[
  x_{0},x_{1},x_{2},\;
  x_{11},x_{12},x_{13},\;
  x_{21},x_{23},x_{32}
\bigr], \]
 with \[\deg(x_{0},x_{1},x_{2},x_{11},x_{12},x_{13},x_{21},x_{23},x_{32})
        =(2,3,3,3,4,5,5,7,11),
\]
and let \(P=\operatorname{Proj}R\cong\mathbb P(2,3,3,3,4,5,5,7,11)\).

Then the following sparse representation of equations gives the quasismoothness of this family.
\[
M \;=\;
\begin{pmatrix}
  x_{11} & x_{12} & x_{13}\\[2pt]
  x_{21} & f_{6}  & x_{23}\\[2pt]
  f_{10} & x_{32} & f_{12}
\end{pmatrix},
\qquad
X \;=\; V\!\bigl(\text{all }2\times2\text{ minors of }M\bigr)\subset P.
\]

where we can choose
\[
\begin{array}{lcl}
f_{6}  &=&
        2x_{0}^{3}
      + x_{1}^{2}
      + 2x_{2}^{2}
      + x_{11}^{2}
      + 3x_{0}x_{12},
\\[6pt]
f_{10} &=&
        x_{0}^{5}
      + 3x_{0}^{2}x_{1}^{2}
      + 3x_{0}^{2}x_{2}^{2}
      + 3x_{0}^{2}x_{11}^{2}
      + 3x_{1}^{2}x_{12}\\&&
      + 2x_{2}^{2}x_{12}
      + x_{11}^{2}x_{12}
      + 2x_{0}x_{12}^{2}
      + 2x_{13}^{2}\\&&
      + 2x_{21}^{2}
      + 2x_{1}x_{23}
      + 2x_{2}x_{23}
      + 2x_{11}x_{23}, \textrm{ and }
\\[6pt]
f_{12} &=&
        2x_{0}^{6}
      + x_{1}^{4}
      + x_{2}^{4}
      + 3x_{11}^{4}
      + 2x_{12}^{3}
      + x_{1}x_{12}x_{13}
      + x_{2}x_{12}x_{13}\\&&
      + 2x_{11}x_{12}x_{13}
      + x_{0}x_{13}^{2}
      + x_{1}x_{12}x_{21}
      + 2x_{2}x_{12}x_{21}
      + x_{11}x_{12}x_{21}\\&&
      + 2x_{0}x_{21}^{2}
      + 3x_{0}x_{1}x_{23}
      + 2x_{0}x_{2}x_{23}
      + 3x_{0}x_{11}x_{23}\\&&
      + 2x_{13}x_{23}
      + 2x_{21}x_{23}.
\end{array}
\]

% \subsection{Terminal (isolated) Fano 4-folds with $|{-}K_X|$ a single elephant}
% \label{Fano4-K-one}
% The following table compiles data   of terminal Fano 4-folds with a single elephant. 
%\begin{table}[h]
% \caption{Terminal Fano 4-folds with single elephant}
% \begin{tabular}{|c|c|c|}\hline
% Format-q_\max & \#Candidates Fano 4 & #\#QS Fano4 with no CY3\ \\\hline
% CI2 & 351 & ? \\\hline
% CI3 & 43 & ?\\\hline
% CI4 & 0 & 0 \\\hline
% $\Gr(2,5)$ & 77 & 24 \\\hline
% $\PxP$ & 18 & 7? \\\hline
% 
% \end{tabular}
% \end{table}

\section{Isolated terminal Fano 4-folds  with  \(h^0(-K_X) \ge 2\)}
\label{Sec:Fano4-NoCY3}
This section consists of the sample examples of ITF4s obtained from  Theorem
\ref{Th:type-II}.  We  present examples that represent the extreme cases, i.e.  those with maximum \(m\) such that \(h^{0}(-K_X)=m \) but a linear section is not a cCY3. We also summarize the the results obtained from the computer search algorithm using \cite{QJSC,BKZ} and compare it to \cite{BKZ}. 
\subsection{Summary}  In this section we summarize the results where \(h^0(-K_X)\ge 2\) but generic  linear section is not an isolated cCY3  in each Gorenstein format, in Table \ref{Combined-Summary-Fano4-Sorted}.  The column labeled ``\# Fano 4-folds
$X^4_{ca}$'' denotes the total number of candidate ITF4s given value of \(h^0(-K_{X})\), up to the sum of the weights  \(W\) that is listed in column ``Format''. The ``\#QS
CY3'' column provides the number of terminal ITF4s such
that  it contains a  linear section which is one of the quasismooth canonical Calabi-Yau 3-folds  obtained by  Brown--Kasprzyk--Zhou in \cite{BKZ}. The column ``\#Fano 4-folds $X^4_{ca}$ NcCY3'' refers to the number of candidate ITF4s  for which a general linear section can not be a  isolated canonical Calabi-Yau 3-fold (cCY3). Finally, the column ``\#QS-$K_2$'' represents the number of quasismooth Fano 4-folds
of Type-$K_2$,  constructed in this article. The case of \(\PxP\) is not dealt
in \cite{BKZ} so one of the columns is empty, though a small number of Calabi--Yau
3-fold  has been constructed in \cite{MNQ-JAL}, that are linear sections of   \(\PxP\) families of ITF4s we found by using our computer search.

\begin{table*}[htb]
\centering
\renewcommand*{\arraystretch}{1.3}
\centering
\caption{Combined summary of Fano 4-folds such that \(h^0(-K_X) \ge 2\) and
the linear section is not an isolated canonical Calabi--Yau 3-fold}
\label{Combined-Summary-Fano4-Sorted}
\begin{tabular}{cccccc}
\toprule
Format & $h^0(-K_X)$ & \#Fano $X^4_{ca}$ & \#QS CY3\cite{BKZ}  & \#Fano $X^4_{ca}$ NcCY3 & \#QS-$K_2$\\ 
\midrule
\evnrow C.I cod.2 
& 2 & 198 & 180 & 13 & 12 \\ 
\(W\)=101& \(3,\ldots,7\) & 73 & 71 & 0 & 0 \\ 
\midrule
\evnrow C.I cod.3  
& 2 & 22 & 9 & 13 & 7\\ 
\(W\)=70& $3,\ldots,8$ & 5 & 4 & 1 & 1 \\ 
\midrule
\evnrow C.I cod.4  
& 2 & 6 & - & 6 & 0\\ 
\(W\)=65& $3$ & 1 & - & 1 & 0 \\ 
\midrule
\evnrow \(\Gr(2,5)\)  
& 2 & 104 & 64 & 20 & 6\ \\ 
\(W=\)70& 3 & 75 & 65 & 6 & 1\\ 
\midrule
\evnrow\(\PxP\)  
& 2 & 58 & - & 30 & 3\\ 
\(W=\)57& 3 & 56 & - & 10 & 1\ \\ 
\evnrow& 4 & 21 & - & 4 & 1 \\ 
\bottomrule
\end{tabular}
\end{table*}

\begin{table}[h]

\end{table}

%%%%%%%%%%%%%%%%%%%%%%%%%%%
\subsection{Sample Examples}
The following are two examples with maximal value of \(h^0(-K_X)\) in both \(\Gr(2,5)\) and \(\PxP\) formats. In  first case $h^0(-K_X)=3$ and $h^0(-K_X)=4$ in the  latter. 
\begin{ej}
Let $X_{ 22, 22, 22, 12, 12}\into \PP(1^3,3,7,11^{3})$ be a regular pullback with weight matrix  
$$\left(\begin{matrix} 1&1&11&11\\ &1&11&11\\ &&11&11\\&&&21 \end{matrix}\right).$$ Then $X$ is a wellformed and quasismooth terminal Fano 4-fold with the following basket of singularities.

$\sB=\{1/3(1, 2, 2, 2),  1/7(3,4,4,4), 3\times 1/11(1,1,3,7)\}.$
 \end{ej}
\begin{ej}
 For a choice of input parameter \(w=(0,0,10,1,1,11)\), a terminal Fano 4-fold \(X\into\PP(1^4,3,7,11^3)\) is given by 
 the weight matrix   $$\left(\begin{matrix} 1&1&11\\ 1&1&11\\ 11&11&21 \end{matrix}\right).$$
Clearly, \(X\) contains singular points of type \[\frac 13(1,2,2,2), \frac 17(3,4,4,4).\] 
In addition, on the weight 11 singular locus, we get 2 points of type \(1/11(1,1,3,7)\). Thus $X$ contains the following basket of singularities,
$$\sB=\{1/3(1, 2, 2, 2),  1/7(3,4,4,4), 2 \times 1/11(1,1,3,7)\}.$$
\end{ej}

%%-----------------------------------------------------

\section*{Acknowledgements} I am grateful   to Gavin Brown and Alexander
Kasprzyk for helpful discussions. 

\bibliographystyle{amsalpha}
\bibliography{References}

\newcommand{\etalchar}[1]{$^{#1}$}
\providecommand{\bysame}{\leavevmode\hbox to3em{\hrulefill}\thinspace}
\providecommand{\MR}{\relax\ifhmode\unskip\space\fi MR }
\providecommand{\MRhref}[2]{%
  \href{http://www.ams.org/mathscinet-getitem?mr=#1}{#2}
}
\providecommand{\href}[2]{#2}
\begin{thebibliography}{KMM{\etalchar{+}}92}

\bibitem[ABR02]{ABR}
S.~Alt{\i}nok, G.~D. Brown, and M.~A. Reid, \emph{Fano 3-folds, {$K3$} surfaces
  and graded rings}, Topology and geometry: commemorating {SISTAG}, Contemp.
  Math., vol. 314, Amer. Math. Soc., Providence, RI, 2002, pp.~25--53.

\bibitem[Bat99]{Bat-toric-4}
Victor~V Batyrev, \emph{On the classification of toric {F}ano 4-folds}, Journal
  of Mathematical Sciences \textbf{94} (1999), no.~1, 1021--1050.

\bibitem[BCP97]{magma}
Wieb Bosma, John Cannon, and Catherine Playoust, \emph{The {M}agma algebra
  system. {I}. {T}he user language}, J. Symbolic Comput. \textbf{24} (1997),
  no.~3-4, 235--265, Computational algebra and number theory (London, 1993).

\bibitem[BK]{GRDB}
G.~Brown and A.~M. Kasprzyk, \emph{Graded ring database, available at},
  http://grdb.co.uk/.

\bibitem[BK16]{BK-ths4}
Gavin Brown and Alexander Kasprzyk, \emph{Four-dimensional projective orbifold
  hypersurfaces}, Experimental Mathematics \textbf{25} (2016), no.~2, 176--193.

\bibitem[BK22]{GRDB-paper}
Gavin Brown and Alexander~M Kasprzyk, \emph{Kawamata boundedness for fano
  threefolds and the graded ring database}, arXiv preprint arXiv:2201.12345
  (2022).

\bibitem[BKQ18]{BKQ}
G.~Brown, A.~M. Kasprzyk, and M.~I. Qureshi, \emph{Fano 3-folds in
  $\mathbb{P}^2\times \mathbb{P}^2$ format, {Tom} and {J}erry}, European
  Journal of Mathematics \textbf{4} (2018), no.~1, 51--72.

\bibitem[BKR12]{BKR}
G.D. Brown, M.~Kerber, and M.A. Reid, \emph{{Fano 3-folds in codimension 4, Tom
  and Jerry. I.}}, Compos. Math. \textbf{148} (2012), no.~4, 1171--1194.

\bibitem[BKZ22]{BKZ}
Gavin Brown, Alexander~M. Kasprzyk, and Lei Zhu, \emph{Gorenstein formats,
  canonical and {C}alabi--{Y}au threefolds}, Experimental Mathematics
  \textbf{31} (2022), no.~1, 146--164.

\bibitem[BRZ13]{BRZ}
A.~Buckley, M.~Reid, and S.~Zhou, \emph{Ice cream and orbifold
  {R}iemann--{R}och}, Izvestiya: Mathematics \textbf{77:3} (2013), 461--486.

\bibitem[CD20]{CD-Fano}
Stephen Coughlan and Tom Ducat, \emph{Constructing fano 3-folds from cluster
  varieties of rank 2}, Compositio Mathematica \textbf{156} (2020), no.~9,
  1873–1914.

\bibitem[CGKS18]{Q-periods}
Tom Coates, Sergey Galkin, Alexander Kasprzyk, and Andrew Strangeway,
  \emph{Quantum periods for certain four-dimensional {F}ano manifolds},
  Experimental Mathematics (2018), 1--39.

\bibitem[CKP15]{CKP}
T.~Coates, A.~Kasprzyk, and T.~Prince, \emph{Four-dimensional {F}ano toric
  complete intersections}, Proceedings of the Royal Society Series A
  \textbf{471} (2015), no.~2175, 20140704, 14. \MR{3303391}

\bibitem[CR02]{wg}
A.~Corti and M.~Reid, \emph{Weighted {G}rassmannians}, Algebraic geometry
  (M.~C. Beltrametti, F.~Catanese, C.~Ciliberto, A.~Lanteri, and C.~Pedrini,
  eds.), de Gruyter, Berlin, 2002, pp.~141--163.

\bibitem[Fuj80]{Fuj1}
Takao Fujita, \emph{On the structure of polarized manifolds with total
  deficiency one, {I}}, Journal of the Mathematical Society of Japan
  \textbf{32} (1980), no.~4, 709--725.

\bibitem[Fuj90]{Fuj4}
\bysame, \emph{Classification theories of polarized varieties}, London Math.
  Soc. Lecture Note Ser. \textbf{155} (1990).

\bibitem[HLM19]{HLM}
Juergen Hausen, Antonio Laface, and Christian Mauz, \emph{On smooth {F}ano
  fourfolds of {P}icard number two}, arXiv preprint arXiv:1907.08000 (2019).

\bibitem[IF00]{fletcher}
A.~R. Iano-Fletcher, \emph{Working with weighted complete intersections},
  Explicit Birational Geometry of 3-folds, vol. 281, London Math. Soc. Lecture
  Note Ser, CUP, 2000, pp.~101--173.

\bibitem[IP99]{Isk-Pro}
V.~A. Iskovskikh and Yu.~G. Prokhorov, \emph{{F}ano varieties}, Algebraic
  geometry, {V}, Encyclopaedia Math. Sci., vol.~47, Springer, Berlin, 1999,
  pp.~1--247. \MR{1668579}

\bibitem[Isk77]{Isk1}
VA~Iskovskih, \emph{Fano 3-folds. {I}}, Mathematics of the USSR-Izvestiya
  \textbf{11} (1977), no.~3, 485.

\bibitem[Isk78]{Isk2}
\bysame, \emph{Fano 3-folds. {II}}, Mathematics of the USSR-Izvestiya
  \textbf{12} (1978), no.~3, 469.

\bibitem[Kal19]{Klash}
Elana Kalashnikov, \emph{Four-dimensional {F}ano quiver flag zero loci},
  Proceedings of the Royal Society A \textbf{475} (2019), no.~2225, 20180791.

\bibitem[Kas13]{Al-twps4}
Alexander Kasprzyk, \emph{Classifying terminal weighted projective space},
  arXiv preprint arXiv:1304.3029 (2013).

\bibitem[KMM{\etalchar{+}}92]{KMM}
J{\'a}nos Koll{\'a}r, Yoichi Miyaoka, Shigefumi Mori, et~al., \emph{Rational
  connectedness and boundedness of {F}ano manifolds}, Journal of Differential
  Geometry \textbf{36} (1992), no.~3, 765--779.

\bibitem[KO{\etalchar{+}}73]{Kob1}
Shoshichi Kobayashi, Takushiro Ochiai, et~al., \emph{Characterizations of
  complex projective spaces and hyperquadrics}, Journal of Mathematics of Kyoto
  University \textbf{13} (1973), no.~1, 31--47.

\bibitem[K{\"u}c95]{Kuchle-95}
Oliver K{\"u}chle, \emph{On {F}ano 4-folds of index 1 and homogeneous vector
  bundles over {G}rassmannians}, Mathematische Zeitschrift \textbf{218} (1995),
  no.~1, 563--575.

\bibitem[K{\" u}c97]{Kuchle-97}
Oliver K{\" u}chle, \emph{Some remarks and problems concerning the geography of
  {F}ano 4-folds of index and {P}icard number one}, Quaestiones Mathematicae
  \textbf{20} (1997), no.~1, 45--60.

\bibitem[MM82]{MM}
S.~Mori and S.~Mukai, \emph{Classification of {F}ano {$3$}-folds with
  {$B_{2}\geq 2$}}, Manuscripta Math. \textbf{36} (1981/82), no.~2, 147--162.

\bibitem[MNQ24]{MNQ-JAL}
Sumayya Mohsin, Shaheen Nazir, and Muhammad~Imran Qureshi, \emph{Construction
  and deformations of {C}alabi--{Y}au 3-folds in codimension 4}, Journal of
  Algebra \textbf{657} (2024), 773--803.

\bibitem[Qur17]{QJSC}
Muhammad~Imran Qureshi, \emph{Computing isolated orbifolds in weighted flag
  varieties}, Journal of Symbolic Computation \textbf{79, Part 2} (2017), 457
  -- 474.

\bibitem[Qur19]{QMOC}
Muhammad~Imran Qureshi, \emph{Biregular models of log del {P}ezzo surfaces with
  rigid singularities}, Mathematics of Computation \textbf{88} (2019), no.~319,
  2497--2521.

\bibitem[Qur21]{QBAM}
\bysame, \emph{Smooth {F}ano 4-folds in {G}orenstein formats}, Bulletin of the
  Australian Mathematical Society \textbf{104} (2021), no.~3, 424--433.

\bibitem[Rei80]{c3f}
M.~Reid, \emph{Canonical {$3$}-folds}, Journ\'ees de {G}\'eometrie
  {A}lg\'ebrique d'{A}ngers, {J}uillet 1979/{A}lgebraic {G}eometry, {A}ngers,
  1979, Sijthoff \& Noordhoff, Alphen aan den Rijn, 1980, pp.~273--310.
  \MR{605348 (82i:14025)}

\bibitem[Rei83]{Reid-MM}
Miles Reid, \emph{Minimal models of canonical 3-folds}, Algebraic varieties and
  analytic varieties, Mathematical Society of Japan, 1983, pp.~131--180.

\bibitem[Tak02]{Takagi}
Hiromichi Takagi, \emph{On classification of \( \mathbb{Q} \)-fano 3-folds of
  {G}orenstein index 2. i, ii}, Nagoya Mathematical Journal \textbf{167}
  (2002), 117--155, 157--216.

\bibitem[Wil87]{Wilson}
PMH Wilson, \emph{Fano fourfolds of index greater than one}, J. reine angew.
  Math \textbf{379} (1987), 172--181.

\bibitem[Wis90]{Wis}
Jaroslaw~A Wisniewski, \emph{Fano 4-folds of index 2 with {$b_2 \ge 2$}. a
  contribution to {M}ukai classification}, Bull. Polish Acad. Sci. Math
  \textbf{38} (1990), 173--184.

\end{thebibliography}
\end{document}